\documentclass[10pt]{article}
\usepackage{amsmath}
\usepackage{amsthm}
\usepackage{amssymb}
\usepackage{url}
\usepackage{latexsym}
\usepackage[xcolor=pst]{pstricks}
\usepackage{graphicx, pst-plot, pst-node, pst-text, pst-tree}
\usepackage{titlefoot}
\usepackage[small]{titlesec}
\usepackage{units} 
\usepackage[small,it]{caption}
\usepackage{rotating} 
\usepackage{cancel}

\usepackage{color}
\definecolor{lightgray}{rgb}{0.8, 0.8, 0.8}
\definecolor{darkgray}{rgb}{0.7, 0.7, 0.7}
\definecolor{darkblue}{rgb}{0, 0, .4}

\usepackage[bookmarks]{hyperref}
\hypersetup{
        colorlinks=true,
        linkcolor=darkblue,
        anchorcolor=darkblue,
        citecolor=darkblue,
        urlcolor=darkblue,
        pdfpagemode=UseThumbs,
        pdftitle={Two vignettes on full rook placements},
        pdfsubject={Combinatorics},
        pdfauthor={Bloom and Vatter},
}

\newcounter{todocounter}

\theoremstyle{plain}

\theoremstyle{definition}

\makeatletter
\newfont{\footsc}{cmcsc10 at 8truept}
\newfont{\footbf}{cmbx10 at 8truept}
\newfont{\footrm}{cmr10 at 10truept}
\renewenvironment{abstract}%
                {
                  \begin{list}{}%
                     {\setlength{\rightmargin}{1in}%
                      \setlength{\leftmargin}{1in}}%
                   \item[]\ignorespaces\begin{small}}%
                 {\end{small}\unskip\end{list}}

\newcommand{\Av}{\operatorname{Av}}
\newcommand{\R}{\mathcal{R}}
\renewcommand{\L}{\mathcal{L}}
\newcommand{\A}{\mathcal{A}}
\newcommand{\B}{\mathcal{B}}
\newcommand{\C}{\mathcal{C}}
\newcommand{\push}{\textsf{s}}
\newcommand{\transfer}{\textsf{t}}
\newcommand{\pop}{\textsf{p}}
\newcommand{\OEISlink}[1]{\href{http://www.research.att.com/projects/OEIS?Anum=#1}{#1}}
\newcommand{\OEISref}{\href{http://www.research.att.com/\~njas/sequences/}{OEIS}~\cite{sloane:the-on-line-enc:}}

\def\sdwys #1{\xHyphenate#1$\wholeString}
\def\xHyphenate#1#2\wholeString {\if#1$%
\else\say{\ensuremath{#1}}\hspace{2pt}%
\takeTheRest#2\ofTheString
\fi}
\def\takeTheRest#1\ofTheString\fi
{\fi \xHyphenate#1\wholeString}
\def\say#1{\begin{turn}{-90}\ensuremath{#1}\end{turn}}

\newenvironment{twostacks}
{
	\begin{scriptsize}
	\psset{xunit=0.0355in, yunit=0.0355in, linewidth=0.02in}
	\begin{pspicture}(0,3)(32,20)
	\psline{c-c}(0,15)(10,15)(10,5)(14,5)(14,15)(18,15)(18,5)(22,5)(22,15)(32,15)
	\rput[l](-0.5,12.5){\mbox{output}}
	\rput[r](32,12.5){\mbox{input}}
}
{
	\end{pspicture}
	\end{scriptsize}
}

\newcommand{\fillstack}[4]{%
	\rput[l](-0.5,17.5){\ensuremath{#1}}
	\rput[c](12.1, 10){\begin{sideways}{\sdwys{#2}}\end{sideways}}
	\rput[c](20.1, 10){\begin{sideways}{\sdwys{#3}}\end{sideways}}
	\rput[r](32,17.5){\ensuremath{#4}}
}

\newcommand{\stackinput}{%
	\psline[linecolor=darkgray]{c->}(24, 17.5)(20, 17.5)(20, 14)
}

\newcommand{\stacktransfer}{%
	\psline[linecolor=darkgray]{c->}(20, 14)(20, 17.5)(12, 17.5)(12, 14)
}
\newcommand{\stackoutput}{%
	\psline[linecolor=darkgray]{c->}(12, 14)(12, 17.5)(8, 17.5)
}

\newcommand{\Lthree}{\L^{\le 3}}

\renewcommand{\u}{\textsf{u}}
\renewcommand{\d}{\textsf{d}}

\amssubj{05A05, 05A15}

\newpagestyle{main}[\small]{
        \headrule
        \sethead[\usepage][][]
        {\sc Full Rook Placements}{}{\usepage}}

\title{\sc Two Vignettes On Full Rook Placements}

\author{%
\begin{tabular}{ccc}
Jonathan Bloom
&\rule{5pt}{0pt}&
Vincent Vatter\footnotemark[\value{footnote}]\footnote{Vatter's research was sponsored by the National Security Agency under Grant Number H98230-12-1-0207 and the National Science Foundation under Grant Number DMS-1301692.  The United States Government is authorized to reproduce and distribute reprints not-withstanding any copyright notation herein.}\\[-0.25ex]
\small Department of Mathematics
&&
\small Department of Mathematics\\[-0.5ex]
\small Dartmouth College
&&
\small University of Florida\\[-0.5ex]
\small Hanover, New Hampshire USA
&&
\small Gainesville, Florida USA\\[-1.5ex]
\end{tabular}
}

\titleformat{\section}
        {\large\sc}
        {\thesection.}{1em}{}

\date{}

\begin{document}
\maketitle

\pagestyle{main}

\begin{abstract}
Using bijections between pattern-avoiding permutations and certain full rook placements on Ferrers boards, we give short proofs of two enumerative results. The first is a simplified enumeration of the $3124$, $1234$-avoiding permutations, obtained recently by Callan via a complicated decomposition. The second is a streamlined bijection between $1342$-avoiding permutations and permutations which can be sorted by two increasing stacks in series, originally due to Atkinson, Murphy, and Ru\v{s}kuc.
\end{abstract}

\section{Introduction}


This note concerns bijections between pattern-avoiding permutations and pattern-avoiding full rook placements (frps) on Ferrers boards. The first is quite simple to define: the permutation (or pattern) $\sigma$ is \emph{contained} in the permutation $\pi$ (both thought of in one-line notation) if $\pi$ has a subsequence which is \emph{order isomorphic} to $\sigma$ (that is, has the same pairwise comparisons as). In this case we write $\sigma\le\pi$. If $\sigma\not\le\pi$, we say that $\pi$ \emph{avoids} $\sigma$. Given a set $B$ of permutations we denote by $\Av(B)$ the \emph{class} of all permutations which avoid every permutation in $B$ and refer to
$$
\sum_{\pi\in\Av(B)} x^{|\pi|}
$$
as the generating function of $\Av(B)$; here $|\pi|$ denotes the length of $\pi$. By constructing a bijection between $\Av(1342)$ and ``$\beta(0,1)$ trees'', B\'ona~\cite{bona:exact-enumerati:} established that the generating function of $\Av(1342)$ is
\begin{equation}
\label{eqn-gf}
\frac{8x^2+12x-1+(1-8x)^{3/2}}{32x}.\tag{$\dagger$}
\end{equation}
Recently, the first author and Elizalde~\cite{bloom:pattern-avoidan:} reproved this result by constructing a much simpler bijection between a symmetry of this class, $\Av(3124)$, and certain frps. To define the latter takes a bit of preparation.

\begin{figure}[t]
\begin{center}
\begin{tabular}{ccccccc}

	\psset{xunit=0.015in, yunit=0.015in, runit=0.04in}
	\psset{linewidth=0.01in}
	\begin{pspicture}(0,0)(50,60)
	\multips{0}(0,0)(10,0){6}{%
		\psline{c-c}(0,0)(0,50)
	}
	\multips{0}(0,0)(0,10){6}{%
		\psline{c-c}(0,0)(50,0)
	}
	\rput[c](5,35){$\times$}
	\rput[c](15,15){$\times$}
	\rput[c](25,45){$\times$}
	\rput[c](35,5){$\times$}
	\rput[c](45,25){$\times$}
	\end{pspicture}

&
\psset{xunit=0.015in, yunit=0.015in}
\begin{pspicture}(0,0)(10,60)
\rput[c](5,25){$\longrightarrow$}
\uput[90](5,25){$\chi$}
\end{pspicture}
&
 
	\psset{xunit=0.015in, yunit=0.015in, runit=0.04in}
	\psset{linewidth=0.01in}
	\begin{pspicture}(0,0)(50,60)
	\multips{0}(0,0)(10,0){4}{%
		\psline{c-c}(0,0)(0,50)
	}
	\multips{0}(40,0)(10,0){2}{%
		\psline{c-c}(0,0)(0,30)
	}
	\multips{0}(0,0)(0,10){4}{%
		\psline{c-c}(0,0)(50,0)
	}
	\multips{0}(0,40)(0,10){2}{%
		\psline{c-c}(0,0)(30,0)
	}
	\rput[c](5,35){$\times$}
	\rput[c](15,15){$\times$}
	\rput[c](25,45){$\times$}
	\rput[c](35,5){$\times$}
	\rput[c](45,25){$\times$}
	\end{pspicture}

&
\psset{xunit=0.015in, yunit=0.015in}
\begin{pspicture}(0,0)(10,60)
\rput[c](5,25){$\longrightarrow$}
\end{pspicture}
&
 
	\psset{xunit=0.015in, yunit=0.015in, runit=0.04in}
	\psset{linewidth=0.01in}
	\begin{pspicture}(0,0)(58,60)
	\multips{0}(0,0)(10,0){4}{%
		\psline{c-c}(0,0)(0,50)
	}
	\multips{0}(40,0)(10,0){2}{%
		\psline{c-c}(0,0)(0,30)
	}
	\multips{0}(0,0)(0,10){4}{%
		\psline{c-c}(0,0)(50,0)
	}
	\multips{0}(0,40)(0,10){2}{%
		\psline{c-c}(0,0)(30,0)
	}
	\psline[linewidth=0.03in]{c-c}(50,0)(50,30)(30,30)(30,50)(0,50)
	\pscircle*(0,50){1}
	\uput[45](0,50){$0$}
	\pscircle*(10,50){1}
	\uput[45](10,50){$1$}
	\pscircle*(20,50){1}
	\uput[45](20,50){$1$}
	\pscircle*(30,50){1}
	\uput[45](30,50){$2$}
	\pscircle*(30,40){1}
	\uput[45](30,40){$1$}
	\pscircle*(30,30){1}
	\uput[45](30,30){$1$}
	\pscircle*(40,30){1}
	\uput[45](40,30){$1$}
	\pscircle*(50,30){1}
	\uput[45](50,30){$2$}
	\pscircle*(50,20){1}
	\uput[45](50,20){$1$}
	\pscircle*(50,10){1}
	\uput[45](50,10){$1$}
	\pscircle*(50,0){1}
	\uput[45](50,0){$0$}
	\rput[c](5,35){$\times$}
	\rput[c](15,15){$\times$}
	\rput[c](25,45){$\times$}
	\rput[c](35,5){$\times$}
	\rput[c](45,25){$\times$}
	\end{pspicture}

&
\psset{xunit=0.015in, yunit=0.015in}
\begin{pspicture}(0,0)(10,60)
\rput[c](5,25){$\longrightarrow$}
\uput[90](5,25){$\Pi$}
\end{pspicture}
&
 
	\psset{xunit=0.015in, yunit=0.015in, runit=0.04in}
	\psset{linewidth=0.01in}
	\begin{pspicture}(0,-10)(100,30)
	\psline[linewidth=0.03in]{c-c}(0,0)(30,30)(50,10)(70,30)(100,0)
	\pscircle*(0,0){1}
	\uput[90](0,0){$0$}
	\pscircle*(10,10){1}
	\uput[90](10,10){$1$}
	\pscircle*(20,20){1}
	\uput[90](20,20){$1$}
	\pscircle*(30,30){1}
	\uput[90](30,30){$2$}
	\pscircle*(40,20){1}
	\uput[90](40,20){$1$}
	\pscircle*(50,10){1}
	\uput[90](50,10){$1$}
	\pscircle*(60,20){1}
	\uput[90](60,20){$1$}
	\pscircle*(70,30){1}
	\uput[90](70,30){$2$}
	\pscircle*(80,20){1}
	\uput[90](80,20){$1$}
	\pscircle*(90,10){1}
	\uput[90](90,10){$1$}
	\pscircle*(100,0){1}
	\uput[90](100,0){$0$}
	\end{pspicture}

\end{tabular}
\end{center}
\caption{An example of the bijections $\chi$ and $\Pi$.}
\label{fig-frps-example}
\end{figure}
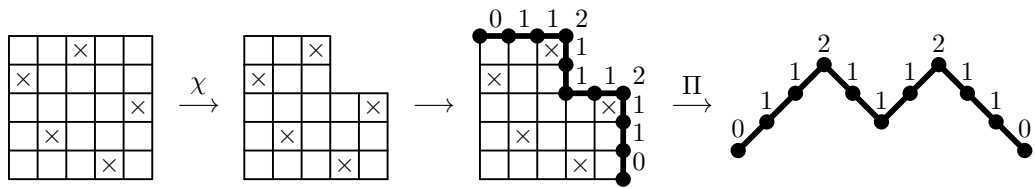

We use French/Cartesian indexing throughout, so for us, a Ferrers board is a left-justified array of cells in which the number of cells in each row is at least the number of cells in the row above. A \emph{full rook placement (frp)} on a Ferrers board consists of a Ferrers board with a designated set of cells, called \emph{rooks}, so that each row and column contains precisely one rook.  For example, the left-most two objects in Figure~\ref{fig-frps-example} are both frps with rooks marked by crosses.

There is a natural partial order on the set of all frps: given two frps $R$ and $S$, we say that $R$ is contained in $S$ if $R$ can be obtained from $S$ by deleting rows and columns. Furthermore, this partial order generalizes the permutation containment order. To make this precise, we call an frp \emph{square} if the underlying Ferrers board is square. There is a natural correspondence between permutations and square frps ($\pi$ maps to the square frp with rooks in the cells $(i,\pi(i))$ for every $i$). When restricted to square frps, the partial order on frps is equivalent to the permutation containment order on the corresponding permutations.

Because of this correspondence between permutations and square frps, we say that an frp \emph{avoids} the permutation $\sigma$ if it avoids the square frp corresponding to $\sigma$. We denote by $\R(B)$ the set of frps that avoid every permutation in $B$ (herein we are only interested in $B=\{312\}$). If we need to stratify this set by number of rooks, we use a subscript, so $\R_n(B)$ consists of the $B$-avoiding frps with precisely $n$ rooks.

A frp is \emph{board minimal} if its rooks do not lie in any smaller Ferrers board, or equivalently, if it has a rook in each of its upper-right corners. We denote by $\R^\times(B)$ the set of $B$-avoiding board minimal frps. Clearly there is a bijection, which we denote by $\chi$, between permutations and board minimal frps. As observed in \cite{bloom:pattern-avoidan:}, because $3124$ ends with its greatest entry and its second-to-last entry is not its second greatest entry, $\chi$ restricts to a bijection from $\Av(3124)$ to $\R^\times(312)$.

In \cite{bloom:a-simple-biject:}, the first author and Saracino constructed a bijection, which we call $\Pi$, from $\R^\times(312)$  to a certain set of labeled Dyck paths denoted by $\L^\times(312)$. The combination of $\chi$ and $\Pi$ is illustrated in Figure~\ref{fig-frps-example}. To construct $\Pi(R)$, we first (as indicated in the third object of Figure~\ref{fig-frps-example}) label every vertex on the northwest-southeast border of the Ferrers board by the length of the longest increasing sequence of rooks lying below and to the left of that vertex. Finally, we rotate this border together with its labels to form a labeled Dyck path.

As shown in \cite{bloom:pattern-avoidan:,bloom:a-simple-biject:}, the labeled Dyck paths in $\L^\times(312)$ are completely characterized by four properties. To state these we need a final definition. A \emph{weak tunnel} in a Dyck path is a horizontal segment between two vertices of the path which stays at or below the path. (As originally defined by Elizalde~\cite{elizalde:fixed-points-an:}, a \emph{tunnel} may only intersect the path at its endpoints, but we have relaxed this condition.) The properties which characterize $\L^\times(312)$ are then:

\begin{itemize}
\item \emph{Monotone property:} labels increase by at most $1$ after an up step and decrease by at most $1$ after a down step.
\item \emph{Zero property:} the $0$ labels are precisely those on the $x$-axis.
\item \emph{Tunnel property:} given two vertices at the same height connected by a weak tunnel, the label of the rightmost vertex is at most the label of the leftmost vertex. (This was called the ``diagonal property'' in \cite{bloom:pattern-avoidan:}.)
\item \emph{Peak property:} the labels rounding any peak (an up step followed by a down step) are $\ell$, $\ell+1$, $\ell$.
\end{itemize}
While not immediately obvious, the mapping $\Pi:R^\times(312)\to \L^\times(312)$ is shown in~\cite{bloom:pattern-avoidan:,bloom:a-simple-biject:} to be bijective.  Thus the problem of enumerating $\Av(3124)$, and thus also $\Av(1342)$, is reduced to the problem of enumerating $\L^\times(312)$. Bloom and Elizalde~\cite{bloom:pattern-avoidan:} were able to show that \eqref{eqn-gf} is indeed the generating function of $\L^\times(312)$.

\section{$\Av(3124, 1234)$}
\label{sec-av-3124-1234}

Permutation classes defined by avoiding two permutations of length four (the ``$2\times 4$ classes'') have proved to be an interesting test bed for comparing enumerative techniques. There are $56$ such classes up to symmetry, and it has been shown that these $56$ classes have $38$ different generating functions (thus some pairs of classes which are inequivalent by symmetry nevertheless lead to the same generating functions, a phenomenon known as \emph{Wilf-equivalence}). Over the past fifteen or so years, authors employing a variety of different approaches  have found all but about a dozen of these generating functions (see Wikipedia~\cite{wikipedia:enumerations-of:} for an up-to-date account).

One of the most recent additions to this catalog is the enumeration of $\Av(3241,4321)$. Verifying a conjecture of Kotesovec, Callan~\cite{callan:permutations-av:} showed that the generating function of this class is
\[
\frac{1}{1-xC(xC(x))},
\]
where $C(x)$ is the generating function for the Catalan numbers (sequence \OEISlink{A165543} in the \OEISref). While Callan's proof is bijective, he concludes his paper by writing ``the bijection presented above works but is hardly intuitive'' and asking ``is there a better proof?'' We provide a positive answer to Callan's question by using frps to enumerate a symmetry of this class, $\Av(3124, 1234)$.  To do so, we use the mapping $\Pi\circ\chi$.  

For concreteness, let us denote an element of $\L^\times_n(312)$ by $(D,\alpha)$, where $D$ is a Dyck path with semilength $n$ (which we denote by $|D|$) and $\alpha$ is the sequence of labels along $D$, ordered from left to right as $\alpha_0$, $\dots$, $\alpha_{2n}$. It is easy to see that the image of $\Av(3124,1234)$ under the mapping $\Pi\circ\chi$ consists of those labeled Dyck paths in $\L^\times(312)$ in which every label is at most $3$. Denoting this set by $\Lthree$, we see that
$$
\sum_{\pi\in\Av(3124,1234)} x^{|\pi|} = \sum_{(D,\alpha)\in \Lthree} x^{|D|}.
$$

The form of the desired generating function suggests that we should decompose labeled Dyck paths in $\Lthree$ based on their returns to the $x$-axis. Suppose that $(D,\alpha)\in\Lthree$ and write it as $\u D^{(1)}\d\u D^{(2)}\d\ldots$ where each $D^{(i)}$ is a (possibly empty) Dyck path. Provided $D^{(i)} \neq \emptyset$, we let $\alpha^{(i)}$ denote the sequence of labels along $D^{(i)}$, decremented by $1$.   Otherwise, we set $\alpha^{(1)} = \emptyset$.  It now follows that the labeled Dyck path $(D^{(i)},\alpha^{(i)})$ satisfies the monotone, tunnel, and peak properties, that every (nonempty) label is between $0$ and $2$, inclusive, and that the zero labels include, but are no longer exclusive to, the $x$-axis. Consequently, we define the set $\A$ to consist of all labeled Dyck paths that arise in this fashion. From the decomposition into returns, it now follows that
$$
\sum_{(D,\alpha)\in \Lthree} x^{|D|}
=
\frac{1}{1-xA(x)},
$$
where $A(x) = \displaystyle \sum x^{|D|}$, the sum taken over all $(D,\alpha)\in \A$.  It only remains to show that $A(x) = C(xC(x))$, or equivalently, that $A(x) = 1+A(x)^2C(x)$.  

We begin by fixing $(D,\alpha)\in \A$ so that $D\neq \emptyset$ and decompose it based on its first return to the $x$-axis as $D= \u D^{(1)}\d D^{(2)}$, where $\alpha^{(1)}$ (respectively, $\alpha^{(2)}$) is the sequence of labels along $D^{(1)}$ (respectively,  $D^{(2)}$). (As before, we define $\alpha^{(i)} =\emptyset$ when $D^{(i)} = \emptyset$.)  Certainly, $(D^{(2)},\alpha^{(2)})\in \A$.  To characterize $(D^{(1)},\alpha^{(1)})$, note that if $\alpha^{(1)}_0 = 0$ or $\alpha^{(1)}=\emptyset$ then  $(D^{(1)},\alpha^{(1)})\in\A$.  On the other hand, if we let $\B$ be the set of all such labeled paths $(D^{(1)}, \alpha^{(1)})$ that arise from this decomposition when $|D^{(1)}|>0$ and $\alpha^{(1)}_0 =1$, i.e., $\B$ has the same characterization as $\A$ except here we insist that its leftmost label is 1, then it (trivially) follows that 
\[
A(x)= 1+x(A(x)+B(x))A(x),
\]
where $B(x) = \displaystyle \sum x^{|D|}$, the sum taken over all $(D,\alpha)\in \B$.   To finish the proof it now suffices to show that $B(x) = A(x)(C(x)-1)$.

To this end we construct a bijection $\phi:\A\times\C\to \B$, where $\C$ is the set of all Dyck paths with positive semilength, which are labeled with all 1s except at the peaks where we place 2s.  (For reference, we call such a labeling \emph{trivial}.)   To define $\phi$, fix $((D,\alpha), E)\in \A\times\C$ and let $i$ be the smallest index so that $\alpha_i=1$; consequently, $D_j = \u$ for all $j\leq i$. We now decompose $E$ as $\u E^{(1)}\d E^{(2)}$ and define $\phi((D,\alpha), E))$ to be the Dyck path
\[
D'=E^{(2)} D_1\ldots D_i  \u E^{(1)}\d  D_i D_{i+1}\ldots,
\]
labeled by
\begin{enumerate}
\item[(a)] placing 1s along the segment $D_1\ldots D_i$,
\item[(b)] giving $\u E^{(1)}\d$ and $E^{(2)}$ the trivial labeling, and
\item[(c)] placing the labels $\alpha_i\alpha_{i+1}\ldots$ along the segment  $D_i D_{i+1}\ldots$.
\end{enumerate}
To see that $\phi$ is bijective first note that  $D=\emptyset$ if and only if the rightmost label on $D'$ is $1$.  In this case we may recover $\u E_1\d$ by decomposing $D'$ by returns.  On the other hand, if $D\neq \emptyset$ then we may recover $\u E^{(1)}\d$ by first decomposing $D'$ into Dyck paths $D^{(1)}D^{(2)}$ so that every label on the $x$-axis in $D^{(1)}$ is labeled with a $1$ and every label on the $x$-axis in $D^{(2)}$, except the first, is labeled with a $0$.  It now follows that $D^{(1)} = E^{(2)}$ and $D^{(2)} = D_1\ldots D_i  \u E^{(1)}\d  D_i D_{i+1}\ldots$.  To see that we may recover $\u E^{(1)}\d$, note that $D_1\ldots D_i$ is a sequence of up-steps that are labeled $0,\ldots, 0,1$ in $(D,\alpha)$.  Consequently, $\u E^{(1)}\d$ is the largest trivially labeled sub-Dyck path in $D^{(2)}$ that contains the leftmost peak of $D^{(2)}$, and has only one return.

\section{Sorting with Two Increasing Stacks in Series}
\label{sec-inc-inc-sorting}

A stack is a last-in first-out sorting device with \emph{push} (move the next entry from the input to the top of the stack) and \emph{pop} (move the entry on top of the stack to the output) operations.  In Volume 1 of \emph{The Art of Computer Programming}~\cite[Section 2.2.1]{knuth:the-art-of-comp:1}, Knuth showed that the permutation $\pi$ can be sorted (meaning that by applying push and pop operations to the sequence of entries $\pi(1),\dots,\pi(n)$ one can output the sequence $1,\dots,n$) if and only if $\pi$ avoids $231$.

Following \emph{The Art of Computer Programming}, several authors, including Knuth himself in \emph{Volume 3}~\cite[Section 5.2.4]{knuth:the-art-of-comp:3}, have studied networks of sorting machines. In particular, the machine consisting of two stacks in series has been intensely analyzed, albeit with limited success. This machine\footnote{It should be noted that another, much more restricted, definition of sorting with stacks in series has been given by West~\cite{west:sorting-twice-t:}. Under this definition, the permutations sortable by two stacks in series do not form a permutation class. These permutations were first counted by Zeilberger~\cite{zeilberger:a-proof-of-juli:}.} allows three operations:
\begin{itemize}
\item \emph{push} the next entry from the input to the top of the first stack, denoted by $\push$,
\item \emph{transfer} the top entry on the first stack to the top of the second stack, denoted by $\transfer$, and
\item \emph{pop} the top entry from the second stack to the output, denoted by $\pop$.
\end{itemize}
Even the problem of determining whether a given permutation can be sorted by this machine has proved to be difficult; Pierrot and Rossin~\cite{pierrot:2-stack-sorting:} have only very recently showed that this problem lies in $\textsf{P}$ (the amount of time required to determine the answer is bounded by a polynomial in the length of the permutation). For the enumeration problem only rough bounds are known, the most recent of which are due to Albert, Atkinson, and Linton~\cite{albert:permutations-ge:}.

\begin{figure}[t]
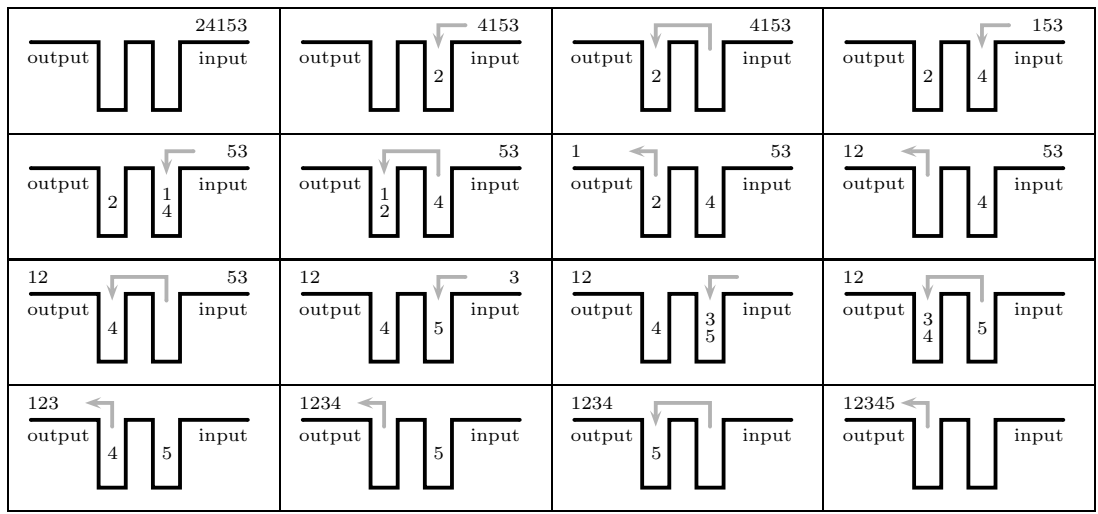

\begin{center}

\begin{tabular}{|c|c|c|c|}
\hline
\begin{twostacks}
\fillstack{}{}{}{24153}
\end{twostacks}
&
\begin{twostacks}
\fillstack{}{}{2}{4153}
\stackinput
\end{twostacks}
&
\begin{twostacks}
\fillstack{}{2}{}{4153}
\stacktransfer
\end{twostacks}
&
\begin{twostacks}
\fillstack{}{2}{4}{153}
\stackinput
\end{twostacks}
\\\hline
\begin{twostacks}
\fillstack{}{2}{41}{53}
\stackinput
\end{twostacks}
&
\begin{twostacks}
\fillstack{}{21}{4}{53}
\stacktransfer
\end{twostacks}
&
\begin{twostacks}
\fillstack{1}{2}{4}{53}
\stackoutput
\end{twostacks}
&
\begin{twostacks}
\fillstack{12}{}{4}{53}
\stackoutput
\end{twostacks}
\\\hline
\begin{twostacks}
\fillstack{12}{4}{}{53}
\stacktransfer
\end{twostacks}
&
\begin{twostacks}
\fillstack{12}{4}{5}{3}
\stackinput
\end{twostacks}
&
\begin{twostacks}
\fillstack{12}{4}{53}{}
\stackinput
\end{twostacks}
&
\begin{twostacks}
\fillstack{12}{43}{5}{}
\stacktransfer
\end{twostacks}
\\\hline
\begin{twostacks}
\fillstack{123}{4}{5}{}
\stackoutput
\end{twostacks}
&
\begin{twostacks}
\fillstack{1234}{}{5}{}
\stackoutput
\end{twostacks}
&
\begin{twostacks}
\fillstack{1234}{5}{}{}
\stacktransfer
\end{twostacks}
&
\begin{twostacks}
\fillstack{12345}{}{}{}
\stackoutput
\end{twostacks}
\\\hline
\end{tabular}
\vspace{0.2in}

\caption{Sorting the permutation $24153$ with two increasing stacks in series.}
\label{fig-incinc-sorting}
\end{center}
\end{figure}

Given the apparent difficulty of analyzing this machine, several authors have considered restricted variants. In particular, Atkinson, Murphy, and Ru\v{s}kuc~\cite{atkinson:sorting-with-tw:} studied sorting with two {\it increasing\/} stacks in series, i.e., two stacks whose entries must be in increasing order when read from top to bottom\footnote{Even without this restriction, the final stack must be increasing from top to bottom if the sorting is to be successful.}. An example of sorting with this machine is shown in Figure~\ref{fig-incinc-sorting}. They proved that this class is characterized by an \emph{infinite} set of minimal avoided permutations, but is nevertheless in bijection with $\Av(1342)$, the most bizarre Wilf-equivalence known to-date. In this section we describe a simple bijection between these permutations and $\L^\times(312)$.

In \cite{atkinson:sorting-with-tw:}, Atkinson, Murphy, and Ru\v{s}kuc associated permutations sortable by two increasing stacks with \emph{greedy stack words}. As our bijection is between $\L^\times(312)$ and these words, we briefly review their definition. First, a \emph{valid stack word} is a word arising from a two-stack sorting (here we do not require that both stacks be increasing --- this condition is imposed later) in which the $i$th letter is $\push$, $\transfer$, or $\pop$, corresponding to whether the $i$th operation is a push, transfer, or pop, respectively. For example, the sorting shown in Figure~\ref{fig-incinc-sorting} corresponds to the valid stack word $\push\transfer\push\push\transfer\pop\pop\transfer\push\push\transfer\pop\pop\transfer\pop$. It is easy to see that the set of \emph{valid} stack words (again, for two not-necessarily increasing stacks in series) are characterized by two rules:
\begin{enumerate}
\item[(W1)] the word contains an equal number of letters equal to $\push$, $\transfer$, and $\pop$, and
\item[(W2)] in every prefix of the word, there are at least as many occurrences of $\push$ as there are of $\transfer$, and at least as many occurrences of $\transfer$ as there are of $\pop$.
\end{enumerate}
To restrict to \emph{increasing} stacks in series, we must impose a further condition:
\begin{enumerate}
\item[(W3)] the word cannot contain a factor (contiguous subsequence) of the form $\transfer u\transfer$ where $u$ is a (possibly empty) valid stack word.
\end{enumerate}
Finally, these rules allow multiple sortings of some permutations, so we imply two \emph{greediness conditions}:
\begin{enumerate}
\item[(W4)] the word cannot contain an $\push\pop$ factor, and
\item[(W5)] the word cannot contain a $u\transfer$ factor where $u$ is a nonempty valid stack word. (The only restriction imposed by (W3) but not by (W5) is that a greedy word representing a sorting by two increasing stacks in series cannot contain a $\transfer\transfer$ factor.)
\end{enumerate}

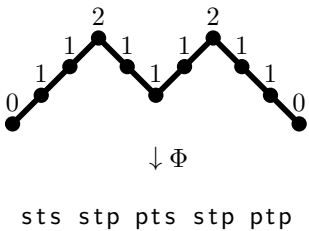
\begin{figure}[t]
\begin{center}

\psset{xunit=0.015in, yunit=0.015in, runit=0.04in}
\psset{linewidth=0.01in}
\begin{pspicture}(0,-35)(100,40)
\psline[linewidth=0.03in]{c-c}(0,0)(30,30)(50,10)(70,30)(100,0)
\pscircle*(0,0){1}
\uput[90](0,0){$0$}
\pscircle*(10,10){1}
\uput[90](10,10){$1$}
\pscircle*(20,20){1}
\uput[90](20,20){$1$}
\pscircle*(30,30){1}
\uput[90](30,30){$2$}
\pscircle*(40,20){1}
\uput[90](40,20){$1$}
\pscircle*(50,10){1}
\uput[90](50,10){$1$}
\pscircle*(60,20){1}
\uput[90](60,20){$1$}
\pscircle*(70,30){1}
\uput[90](70,30){$2$}
\pscircle*(80,20){1}
\uput[90](80,20){$1$}
\pscircle*(90,10){1}
\uput[90](90,10){$1$}
\pscircle*(100,0){1}
\uput[90](100,0){$0$}
\rput[B](5,-35){$\push$}
\rput[B](10,-35){$\transfer$}
\rput[B](15,-35){$\push$}
\rput[B](25,-35){$\push$}
\rput[B](30,-35){$\transfer$}
\rput[B](35,-35){$\pop$}
\rput[B](45,-35){$\pop$}
\rput[B](50,-35){$\transfer$}
\rput[B](55,-35){$\push$}
\rput[B](65,-35){$\push$}
\rput[B](70,-35){$\transfer$}
\rput[B](75,-35){$\pop$}
\rput[B](85,-35){$\pop$}
\rput[B](90,-35){$\transfer$}
\rput[B](95,-35){$\pop$}
\rput(50,-12){$\downarrow$}
\uput[0](50,-12){$\Phi$}
\end{pspicture}

\end{center}
\caption{An example of the bijection between $\L^\times(312)$ and greedy stack words.}
\label{fig-stack-bij-example}
\end{figure}

Our goal, then, is to construct a bijection, which we call $\Phi$, from labeled Dyck paths in $\L^\times(312)$ to stack words satisfying (W1)--(W5). To explain this bijection we introduce a bit of new terminology: a step is \emph{positive} (resp., \emph{neutral}, \emph{negative}) if the label of the vertex it leads to is greater than (resp., equal to, less than) the label of the preceding vertex. The monotone property shows that labels can change by at most $1$ during a step, and that there are only four types of steps: positive and neutral up steps and neutral and negative down steps. We map each of these four types of steps to one or two letters in $\Phi(D,\alpha)$. The four cases are shown below, while a complete example of the bijection is given in Figure~\ref{fig-stack-bij-example}.
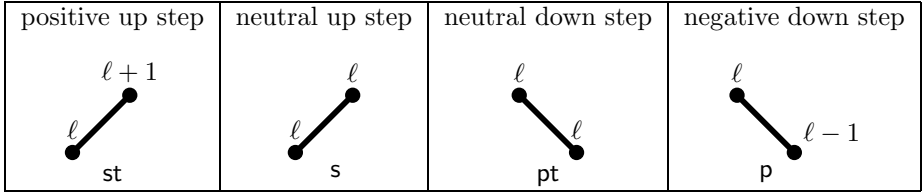
\begin{center}
\begin{tabular}{|c|c|c|c|}
\hline
\mbox{positive up step}
&
\mbox{neutral up step}
&
\mbox{neutral down step}
&
\mbox{negative down step}\\[10pt]
\psset{xunit=0.030in, yunit=0.030in, runit=0.04in}
\psset{linewidth=0.01in}
\begin{pspicture}(0,-5)(14,16)
\psline[linewidth=0.03in]{c-c}(0,0)(10,10)
\pscircle*(0,0){1}
\uput[90](0,0){$\ell$}
\pscircle*(10,10){1}
\uput[90](10,10){$\ell+1$}
\uput[-90](7,0){$\push\transfer$}
\end{pspicture}
&
\psset{xunit=0.030in, yunit=0.030in, runit=0.04in}
\psset{linewidth=0.01in}
\begin{pspicture}(0,-5)(10,16)
\psline[linewidth=0.03in]{c-c}(0,0)(10,10)
\pscircle*(0,0){1}
\uput[90](0,0){$\ell$}
\pscircle*(10,10){1}
\uput[90](10,10){$\ell$}
\uput[-90](7,0){$\push$}
\end{pspicture}
&
\psset{xunit=0.030in, yunit=0.030in, runit=0.04in}
\psset{linewidth=0.01in}
\begin{pspicture}(0,-5)(10,16)
\psline[linewidth=0.03in]{c-c}(0,10)(10,0)
\pscircle*(0,10){1}
\uput[90](0,10){$\ell$}
\pscircle*(10,0){1}
\uput[90](10,0){$\ell$}
\uput[-90](5,0){$\pop\transfer$}
\end{pspicture}
&
\psset{xunit=0.030in, yunit=0.030in, runit=0.04in}
\psset{linewidth=0.01in}
\begin{pspicture}(0,-5)(20,16)
\psline[linewidth=0.03in]{c-c}(0,10)(10,0)
\pscircle*(0,10){1}
\uput[90](0,10){$\ell$}
\pscircle*(10,0){1}
\uput[45](10,0){$\ell-1$}
\uput[-90](5,0){$\pop$}
\end{pspicture}
\\\hline
\end{tabular}
\end{center}
Note that by the construction of $\Phi$ --- assuming $\Phi(D,\alpha)$ corresponds to a valid stack word --- the label of a vertex of $(D,\alpha)$ will equal the number of entries in the second stack at the corresponding point in the sorting.

Let $(D,\alpha)\in\L^\times_n(312)$ be arbitrary. We aim to show that the word $w=\Phi(D,\alpha)$ satisfies (W1)--(W5). First, note that $w$ contains $n$ occurrences of both $\push$ and $\pop$ because $D$ is a Dyck path. Suppose that $(D,\alpha)$ contains $a^+$ positive up steps. By the zero property, it must then contain precisely $a^+$ negative down steps, or equivalently, $n-a^+$ neutral down steps. This shows that $w$ also contains $n$ occurrences of the letter $\transfer$, so $w$ satisfies property (W1).

To check that $\Phi(D,\alpha)$ satisfies (W2), consider an arbitrary prefix $u$ of $w$. This prefix corresponds to an initial segment of $(D,\alpha)$, which we denote by $(E,\beta)$. While not necessarily a Dyck path, $(E,\beta)$ still satisfies the monotone, zero, tunnel, and peak properties. Suppose that $u$ contains $a$ occurrences of $\push$, $b$ occurrences of $\transfer$, and $c$ occurrences of $\pop$, so we want to show that $a\ge b\ge c$. If $(D,\alpha)$ has $a^+$ positive up steps then by the zero property, it contains at most $a^+$ negative down steps, and thus it contains at least $c-a^+$ neutral down steps. This shows that $b\ge a^+ + c-a^+=c$, which is one of the inequalities needed for (W2).

The other inequality we need to show is $a\ge b$. Suppose to the contrary that this inequality fails for some prefix $u$ of $w$ and choose $u$ to be as short as possible subject to this constraint. By the minimality of $u$, it must end in $\pop\transfer$, and thus the last step in $(E,\beta)$ is a neutral down step. Suppose that this final step ends at the vertex $y$ and let $x$ denote the rightmost vertex to the left of $y$ at the same level as $y$. We break $(E,\beta)$ into two pieces at $x$; suppose that there are $a_1$ up steps and $c_1$ down steps in $(E,\beta)$ to the left of $x$ (so $a_1\ge c_1$) and $a_2$ up steps and $c_2$ down steps between $x$ and $y$ (so $a_2=c_2$). Further suppose there are a total of $b_1$ positive up steps and neutral down steps (these are precisely the steps which correspond to a $\transfer$ in $u$) to the left of $x$ (so, by the minimality of $y$, $a_1\ge b_1$) and $b_2$ positive up steps and neutral down steps between $x$ and $y$. Now if there are $a_2^+$ positive up steps between $x$ and $y$, the tunnel property implies that there are at least $a_2^+$ negative down steps between $x$ and $y$. This shows that
\[
b_2\le a_2^+ + (c_2-a_2^+)=a_2^+ + (a_2-a_2^+)=a_2,
\]
so the number of occurrences of $\transfer$ in $u$ is at most $b_1+b_2\le a_1+a_2$, a contradiction.

To see that (W3) holds, first note that by the construction of $\Phi$ it is clearly impossible for $w$ to contain a $\transfer\transfer$ factor. Next we establish (W5), and thus the rest of (W3). Suppose to the contrary that $w$ does contain a $u\transfer$ factor for a valid stack word $u$. Then the letters of $u\transfer$ correspond to a Dyck subpath of $(D,\alpha)$, and thus the leftmost and rightmost vertices of this subpath are be connected by a weak tunnel. Assume that there are $a$ up steps in this subpath (so there are also $a$ down steps in the subpath), and that $a^+$ of these are positive up steps. By the tunnel property, the label of the rightmost vertex is at most the label of the leftmost vertex, so this subpath contains at least $a^+$ negative down steps, and thus at most $a-a^+$ neutral down steps. Thus $u$ contains at most $a$ occurrences of the letter $\transfer$. However, because this subpath corresponds to a $u\transfer$ factor, the final $\transfer$ must correspond to a final neutral down step. But then $u$ can contain at most $a-1$ occurrences of $\transfer$, and thus is not a valid stack word.

The final property, (W4), follows quickly. Were $w$ to contain an $\push\pop$ factor, it could only be the result of an up step followed immediately by a down step, i.e., a peak. However, the peak property states that peaks can only occur as positive up steps followed by negative down steps, and thus peaks correspond to $\push\transfer\pop$ factors, not $\push\pop$ factors.

The inverse of $\Phi$ is easier to describe. To construct $\Phi^{-1}(w)=(D,\alpha)$, we build a Dyck path in which pushes correspond to up steps and pops correspond to down steps. We then label the vertices of this path by the number of entries in the second stack before the next push or pop. (W1) and (W2) ensure that we do indeed obtain a Dyck path. Because of (W3), $w$ cannot contain consecutive occurrences of $\transfer$, so the labels can increase or decrease by at most $1$ at each step. Clearly labels can only decrease on down steps ($w$ cannot contain a $\pop\transfer\transfer$ factor by (W3)), while (W4) ensures that labels can only increase on up steps, thereby verifying that $(D,\alpha)$ satisfies the monotone property. This labeled Dyck path also satisfies the zero property because at the end of the sorting described by $w$, nothing remains in the second stack. The peak property follows because $w$ cannot contain an $\push\pop$ factor by (W4).

It remains only to check the tunnel property. Suppose to the contrary that $(D,\alpha)$ fails the tunnel property, and choose vertices $x$ and $y$ connected by a weak tunnel with $x$ to the left of $y$ such that the label of $y$ is greater than the label of $x$. Subject to these constraints, further choose $x$ and $y$ to be as close to each other as possible. Clearly if the tunnel connecting $x$ and $y$ touches the path in its interior, then there will be a violation of the tunnel property strictly between $x$ and $y$, a contradiction to our choice of these vertices. Because $x$ and $y$ are connected by a weak tunnel, $x$ must be followed by an up step and $y$ must follow a down step. If the vertex before $y$ is connected to $y$ by a negative down step, then it and the vertex after $x$ violate the tunnel property (they are connected by a weak tunnel because the tunnel between $x$ and $y$ does not touch the path in its interior), contradicting our choice of $x$ and $y$. Thus the vertex before $y$ must be connected to $y$ by a neutral down step, which corresponds to a $\pop\transfer$ factor in $w$. However, this shows that $w$ contains a $u\transfer$ factor for a valid stack word $u$, contradicting (W5).

\bibliographystyle{acm}
\bibliography{../../refs}

\end{document}